%% file: main.tex
\documentclass[letterpaper, 10 pt, conference]{ieeeconf}  

\IEEEoverridecommandlockouts                              
\usepackage{amsmath} 
\usepackage{amssymb}  
\usepackage{hyperref}
\usepackage{mathtools}
\usepackage{bm,fixmath}
\usepackage{tikz}
\usepackage{accents}
\usepackage{siunitx}
\usepackage{pgfplots}
\usepackage{mathrsfs}
\usepackage{graphicx}
\pgfplotsset{compat=newest}
\pgfplotsset{plot coordinates/math parser=false}
\newlength\figureheight
\newlength\figurewidth 
\input{./Figures/tikzFigOptions.tex}

\newtheorem{assum}{Assumption}
\newtheorem{rem}{Remark}
\newtheorem{prop}{Proposition}
\newtheorem{definition}{Definition}[section]

\title{\LARGE \bf Anti-windup design for internal model online constrained  optimization}

\author{Umberto Casti$^{1}$ and Sandro Zampieri$^{2}$
\thanks{*This work was not supported by any organization}
\thanks{$^{1}$Umberto Casti is with the Department of Information Engineering (DEI), University of Padova, 35131 Padova, Italy
        (\href{mailto:castiumber@dei.unipd.it}{\texttt{castiumber@dei.unipd.it}}).}%
\thanks{$^{2}$Sandro Zampieri is with the Department of Information Engineering (DEI), University of Padova, 35131 Padova, Italy (\href{mailto:zampi@dei.unipd.it}{\texttt{zampi@dei.unipd.it}}).}%
}

\input{notation}

\begin{document}

\maketitle
\thispagestyle{empty}
\pagestyle{empty}

\begin{abstract}
This paper proposes a novel algorithmic design procedure for online constrained optimization grounded in control-theoretic principles. By integrating the Internal Model Principle (IMP) with an anti-windup compensation mechanism, the proposed Projected-Internal Model Anti-Windup (P-IMAW) gradient descent exploits a partial knowledge of the temporal evolution of the cost function to enhance tracking performance. The algorithm is developed through a structured synthesis procedure: first, a robust controller leveraging the IMP ensures asymptotic convergence in the unconstrained setting. Second, an anti-windup augmentation guarantees stability and performance in the presence of the projection operator needed to satisfy the constraints. The effectiveness of the proposed approach is demonstrated through numerical simulations comparing it against other classical techniques.
\end{abstract}
\section{INTRODUCTION}
Recent advances in modern technology have revealed the need for increasingly sophisticated and high-performance optimization algorithms able to solve non-stationary (time-varying) optimization problems, that are commonly called \textit{online optimization} problems.
These naturally arise in diverse applications, including control~\cite{paternain_realtime_2019}, signal processing~\cite{natali_online_2021}, and machine learning~\cite{chang_distributed_2020}. These scenarios, that are typically characterized by time-varying cost functions and often also by the presence of constraints, require solving a sequence of optimization problems in real time to ensure applicability in dynamic environments~\cite{dallanese_optimization_2020,simonetto_timevarying_2020}.

A common approach to the solution of these problems consists in adapting static optimization techniques, such as projected gradient descent, to the online setting. These so-called \textit{unstructured} methods have been extensively studied~\cite{selvaratnam_numerical_2018,cao_online_2019,zhang_online_2021,lupien_online_2023}, but they typically provide limited tracking performance as they do not exploit any model of the temporal evolution of the problem.

To address this limitation, \textit{structured} algorithms explicitly incorporate information about the time-varying nature of the problem. Notable examples include prediction-correction schemes~\cite{simonetto_dual_2019,bastianello_extrapolation_2023}, interior-point methods~\cite{fazlyab_prediction_2018}, and virtual queue algorithms~\cite{cao_virtual-queue-based_2018}. These methods rely on assumptions about the rate of change in the cost and in the constraints and offer improved performance over unstructured approaches. However, they do not fully leverage more detailed models of the problem dynamics, thereby limiting tracking precision.

In this work, we propose a novel algorithmic framework for \textit{online optimization} problems with constraints, grounded in control-theoretic principles—most notably the IMP and anti-windup compensation. Our approach is based on an online extension of projection-based methods and is capable of handling general classes of constraints, e.g., projections onto convex sets, although for simplicity we focus our analysis on the representative case of nonnegativity constraints.

Our contribution lies at the intersection of optimization and control theory. Prior works have explored control-theoretic designs in both static~\cite{lessard_analysis_2016,scherer_optimization_2023} and online settings~\cite{casti_control_2024,simonetto_optimization_2023,davydov_contracting_2023,simonetto_nonlinear_2024}. For instance,~\cite{shahrampour_online_2017,shahrampour_distributed_2018} analyze online algorithms under linear dynamics of the optimal trajectory;~\cite{simonetto_optimization_2023,simonetto_nonlinear_2024} employ Kalman filtering in stochastic settings;~\cite{casti_stochastic_2024} applies robust $\mathcal{H}_\infty$ techniques to unconstrained problems;~\cite{casti_control_2024} introduces a discrete-time IMP-based algorithm for equality-constrained problems; and~\cite{davydov_contracting_2023} uses contraction analysis to study continuous-time unstructured methods.

In contrast, the algorithmic synthesis presented in this paper is based on a control framework that leverages robust control techniques and the IMP to address the challenges posed by time variation and partial knowledge of the \textit{online optimization} problem. To enforce nonnegativity constraints, we incorporate an anti-windup mechanism that mitigates saturation effects and enables accurate tracking in constrained settings.

\noindent \textbf{Notation.} We denote by $\RNumb$ and $\RNumb_+$ the set of real and nonnegative numbers, and by $\mathbb{R}\left(z\right)$ the set of rational functions in $s$ with real coefficients. Vectors and matrices are denoted by bold letters, e.g., $\x \in \Rn$ and $\A \in \Rnn$, while the $i$-th component of a vector is denoted by a subscript, e.g., $x_i \in \RNumb$. The identity matrix of order $n$ is denoted by $\I_n$. The $2$- of a vector is denoted by $\norm{\cdot}$. In addition, we adopt the entry-wise partial order in $\Rn$, so that $\x \leq \y$ means the inequality holds componentwise. The notation $\A \preceq \B$ (resp. $\A \prec \B$) indicates that $\B - \A$ is positive semi-definite (resp. positive definite). The Kronecker and Hadamard products are denoted by $\otimes$ and $\odot$, respectively. The $\mathcal{L}_2$ space is defined as in~\cite{kothare_unified_1994}, and the corresponding $\mathcal{L}_2$ gain as in~\cite{grimm_linear_2004}. Finally, we denote the Laplace transform with the symbol $\Laplace$.

\section{PROBLEM FORMULATION}
Consider the following continuous-time, time-varying nonnegative optimization problem:   
\begin{equation}\label{eq:problem1}
\xStart = \argmin_{\x \geq 0} \ftx, \quad t \in \RNumb_+,
\end{equation}
where $\ftx$ denotes the time-varying quadratic cost function
\begin{equation}\label{eq:quadratic-cost} 
\ftx = \frac{1}{2}\x^\top \A\x + \bTt\x + \ct,
\end{equation}
where $\A \in \Rnn$ is such that $\A = \A^\top \succ 0$, and $\ct \in \RNumb$, $\bt \in \Rn$ for all $t$.
Following the approach in~\cite{casti_control_2024}, in this paper we extend the results proposed therein to address the case in which we have the nonnegative constraint $\x \geq 0$.


Due to the convexity of problem~\eqref{eq:problem1}, it is possible to prove that it admits a unique minimizer $\xStart$ for all $t \in \RNumb_+$, which is independent of $\ct$. Therefore, $\ct$ can be set to zero without loss of generality.

\begin{rem}
The online optimization problem \eqref{eq:problem1} with $f_t\left(\x\right)$ given in \eqref{eq:quadratic-cost} can be seen as the approximation of an optimization problem with a more general cost function. Indeed, if we take the second-order Taylor expansion of a general time-varying convex cost function $f_t\left(\x\right)$ around its optimal solution $\xStart$ we get
\begin{equation}\label{eq:fapprox}
f_t(\x)\simeq f_t(\xStart)+\frac{1}{2}(\x-\xStart)^T\Delta f_t(\xStart(\x-\xStart)
\end{equation}
in which the linear term is missing by the optimality of $\xStart$ and where $\Delta f_t(\xStart)$ is the Hessian of $f_t$ evaluated in $\xStart$. 
If we assume that this is almost time-invariant and if we define 
$\A\coloneqq \Delta f_t(\xStart)$, 
$b(t)\coloneqq\A\xStart$ and $c(t)\coloneqq f_t(\xStart)$, we get that \eqref{eq:fapprox} is equivalent to
$$f_t(\x)\simeq \frac{1}{2}\x^\top \A\x + \bTt\x + \ct.$$
We can conclude that \eqref{eq:quadratic-cost} is a good approximation of a more general cost function $f_t(\x)$ if the Hessian of this cost function evaluated in its optimal solution $\xStart$ is close to be time-invariant and if we are able to keep $\x$ close to $\xStart$. This sounds similar to what we do when we approximate a nonlinear model by its linearization, which is a reasonable approximation only in case we are able, by the control action, to keep the  state closed to the state around which the system has been linearized.
\end{rem}

As it is common in optimization, we assume that the optimization algorithm does not have exact knowledge of the objective function $\ftx$ in~\eqref{eq:quadratic-cost}, but that it can be informed about the gradient of $\ftx$ evaluated in prescribed values of $\x$. Moreover, we assume the following:

\begin{assum}\label{as:hessianA}
The matrix $\A$ satisfies:
\begin{equation}\label{eq:estA}
\lm\I_n \preceq \A \preceq \lM\I_n,
\end{equation}
where $0 < \lm \leq \lM < +\infty$ and with $\lm$ and $\lM$ assumed to be known.
\end{assum}

\begin{rem}
    Assumption~\ref{as:hessianA} ensures that the cost function $\ftx$ is $\lm$-strongly convex and $\lM$-smooth for any time $t \in \RNumb_+$~\cite{bastianello_internal_2022}. 
\end{rem}

By adopting Assumption~\ref{as:hessianA}, our work aligns with the existing literature. Furthermore, following the ideas of~\cite{bastianello_internal_2022} and merging standard assumptions in control theory with those in optimization, we assume that the time-varying parameter $\bt$ in \eqref{eq:quadratic-cost} is governed by a Linear Time-Invariant (LTI) Multi-Input Multi-Output (MIMO) system, where each component $\bit$ evolves according to the following autonomous system:

\begin{assum}[Model of $\bt$]\label{as:modelbk}
Each entry $\bit$ of the vector $\bt$ evolves according to the continuous-time Single-Input Single-Output (SISO) LTI system:
\begin{equation}\label{eq:linearModels}
\Sigma_i:\begin{cases}
\dxiit &= \F\xiit, \\
\bit &= \H \xiit,
\end{cases}
\end{equation}
where $\F$ is marginally stable, with state dimension $m$. The initial state $\xiB_i(0)\in \Rm$ is assumed to be unknown.
\end{assum}

Under these hypotheses, the evolution of the vector $\bt$ can equivalently be represented by the MIMO LTI system:
\begin{equation}\label{eq:linearModelsEst}
\Sigma_{ext}:\begin{cases}
\dxiextt &= \Fext\xiextt, \\
\bt &= \Hext \xiextt,
\end{cases}
\end{equation}
with
\begin{equation}
\xiextt=\begin{bmatrix} \xiOnet\\ \vdots\\ \xint \end{bmatrix}\in\RNumb^{mn},
\end{equation}
and matrices $\Fext = \I_n \otimes \F$ and $\Hext =\I_n \otimes \H$.

\begin{rem}
    Assumption~\ref{as:modelbk} implies the existence of a marginally stable model governing the evolution of the term $\bt$ in~\eqref{eq:quadratic-cost}. While the case where $\F$ has all eigenvalues with negative real part is less interesting, since in this case $\bt$ asymptotically vanishes for any initial condition, leading to $\xStart$ converging to zero, the more interesting scenario arises when $\F$ has eigenvalues with zero real part. In this case, $\bt$ does not decay over time and can exhibit complex behaviors such as oscillations, linear trends, and quadratic ramps, making the optimization problem~\eqref{eq:problem1} more interesting.
\end{rem}

\begin{rem} One might argue that Assumption~\ref{as:modelbk} suggests identical evolution among all components $\bit$. This interpretation, however, is misleading. The model~\eqref{eq:linearModelsEst} allows for highly diverse behaviors of $\bt$; for instance, distinct components may evolve as linear ramps with different slopes, some may follow sinusoidal ramps, while others exhibit alternative dynamics. 
\end{rem}

\section{ALGORITHM}

We propose here a first-order optimization method, the \textit{Projected-Internal Model Anti-Windup (P-IMAW)} gradient descent, that assumes to access to an oracle providing the gradient on a given point $\xt$:
\begin{equation}\label{eq:gradf}
    \gradfxt = \A\xt + \bt.
\end{equation}
Our algorithm can be seen as a generalization of the projected gradient descent~\cite{parikh_proximal_2014}, tailored to solve problem~\eqref{eq:problem1}. 

The \textit{Online Projected-Gradient Descent} (OP-GD) can be conveniently interpreted as the following continuous time dynamic system~\cite{xing_exponential_2003}
\begin{equation}\label{eq:onlineProjDesc}
\begin{array}{l}
\dqt = -\qt + \proj_{\geq 0}\left(\pt\right)\\
\pt = \qt - \alpha\gradfqt, 
\end{array}
\end{equation}
where $\alpha \in \RNumb_+$ and $\projPlus$ represents the nonnegative orthant projection operator $\projPlus:\Rn\to \RNumb^n_+$.
\begin{rem}
The nonnegativity constraint in problem~\eqref{eq:problem1} determines the form of the projection operator used in the projected algorithm~\eqref{eq:onlineProjDesc}. The reasoning developed here can be extended to more general projection operators~\cite{parikh_proximal_2014}, provided they satisfy certain properties, such as those in Definition~\ref{def:PhiI}. Under suitable assumptions, the proposed approach can therefore be generalized to other types of constraints, e.g. to general box constraints.
\end{rem}

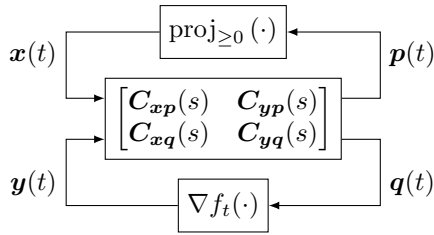
\begin{figure}[ht]
    \vskip 0.0in
    \begin{center}
        \centerline{\IfFileExists{./Figures/controllerGeneral.tex}{\input{./Figures/controllerGeneral.tex}}{?}}
        \caption{
        General online projected gradient descent scheme. 
        }
        \label{fig:algGeneral}
    \end{center}
    \vskip -0.2in
\end{figure}

The standard OP-GD dynamics presented in~\eqref{eq:onlineProjDesc} can be viewed as a special case of the more general online projected gradient descent scheme illustrated in Fig.~\ref{fig:algGeneral}.

Indeed, by setting
\begin{equation*}
\begin{array}{rlrl}
     \Cxps &= \frac{1}{1+s}\In & \Cyps &= 0, \\[0.5em]
     \Cxqs &= \frac{1}{1+s}\In & \Cyqs &= -\alpha\In,
\end{array}
\end{equation*}
the resulting dynamics recover exactly the OP-GD algorithm in~\eqref{eq:onlineProjDesc}.

\begin{figure}[ht] 
\vskip 0.2in 
\begin{center} 
\centerline{\IfFileExists{./Figures/controllerMainSheme.tex}{\input{./Figures/controllerMainSheme.tex}}{?}} 
\caption{\textit{Projected Internal Model Anti-Windup} (P-IMAW) gradient descent scheme proposed in this work. Here, $\Cs$ denotes a robust internal model controller, $\projPlus$ the projection onto the nonnegative orthant, and $\TT$ the anti-windup compensation mechanism. The algorithm assumes access to an oracle providing the gradient $\gradf$.}
\label{fig:algGen} 
\end{center} 
\vskip -0.2in 
\end{figure}
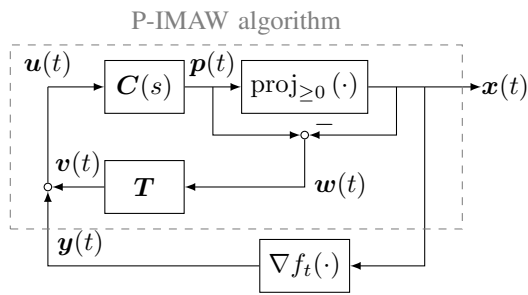

The algorithm proposed here is another instance of the general scheme in Fig.~\ref{fig:algGeneral}. This is obtained by incorporating concepts introduced in~\cite{bastianello_internal_2022}, where an internal model gradient descent algorithm was proposed in the \textit{unconstrained} setting.
Here, we propose a design method extending the one proposed in~\cite{bastianello_internal_2022} and able to deal with the nonlinearity resulting from the projection step. This method, which is described by the block diagram in Fig.~\ref{fig:algGen}, results from the scheme proposed in~\cite{bastianello_internal_2022} in which, besides adding the projection step, there is also an anti-windup compensation able to attenuate the effects of the projection.
Also this scheme is an instance of the general control framework in Fig.~\ref{fig:algGeneral}, in which transfer matrices are
\begin{equation}\label{eq:genTranferMatr}
\begin{aligned}
    \Cxps &= \left[\In - \Cs \TT\right]^{-1} \Cs \TT, \\
    \Cyps &= \left[\In - \Cs \TT\right]^{-1} \Cs, \\
    \Cxqs &= \In, \qquad \Cyqs = 0.
\end{aligned}
\end{equation}

Our method is able to ensure the stability of the resulting nonlinear closed-loop system and to guarantee a certified performance level (formally characterized in Subsection~\ref{subsec:step2}) for problem~\eqref{eq:problem1}. In particular, we derive a  structured approach comprising one analysis step and two synthesis steps. Specifically, the two-step synthesis procedure aligns with the widely adopted unified framework presented in~\cite{kothare_unified_1994}, based on the following paradigm: \textquotedblleft Design the linear controller ignoring control input nonlinearities and then add Anti-Windup Bumpless Transfer (AWBT) compensation to minimize the adverse effects of any control input nonlinearities on closed-loop performance".

Specifically, in \textit{Step 0} we decouple the online optimization problem by leveraging the structure of the objective function and the eigendecomposition $\A = \V\LambdaB\VT$ so to simplify the analysis and the control synthesis. In \textit{Step 1}, we exploit a robust IMP to design a strictly proper controller $\Cs\in\Rnn{\left(s\right)}$, first addressing the \textit{unconstrained} case and, at the same time, dealing the uncertainty of $\A$. Finally, in \textit{Step 2}, we address the \textit{constrained} case by designing the anti-windup static matrix $\TT\in\Rnn$ in Fig.~\ref{fig:algGen}, which compensates for input saturation effects while preserving stability and performance guarantees.
\subsection{Step 0: Decoupling Analysis}

In order to conveniently manage the partial knowledge we assumed of the matrix $\A$ (see Assumption~\ref{as:hessianA}), we restrict the structure of the controller and of the anti-windup mechanism in Fig.~\ref{fig:algGen} to be such that $\Cs = \cs\In$ and $\TT = \rho\In$, where $\cs\in\RNumb{\left(s\right)}$ scalar strictly proper transfer function and $\rho\in\RNumb$. 
Hence the block diagram in Fig.~\ref{fig:algGen} is equivalent to the block diagram in Fig.~\ref{fig:alg0}. This follows from the observation that, in Fig.~\ref{fig:algGen}, $\ut = \yt + \vt$, and thus
\begin{equation}\label{eq:uct}
    \ut =\left( \A\xt + \bt\right) +\vt =\left( \A\xt + \vt\right) + \bt .
\end{equation}
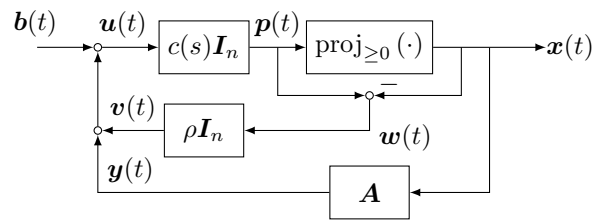
\begin{figure}[ht] 
\vskip 0.2in 
\begin{center} 
\centerline{\IfFileExists{./Figures/controllerTmp.tex}{\input{./Figures/controllerTmp.tex}}{?}} 
\caption{Equivalent block diagram of Fig.~\ref{fig:algGen} after applying the substitutions~\eqref{eq:gradf}, $\Cs = \cs \In$ and $\TT = \rho \In$.} 
\label{fig:alg0} 
\end{center} 
\vskip -0.2in 
\end{figure}
Exploiting the symmetry of the matrix $\A$, we can write $\A = \V\LambdaB\VT$, where $\V$ is a square orthogonal matrix and $\LambdaB$ is a diagonal matrix having the eigenvalues of $A$ in the diagonal. Based on this, the scheme in Fig.~\ref{fig:alg0} can be equivalently rewritten as in Fig.~\ref{fig:alg1}.
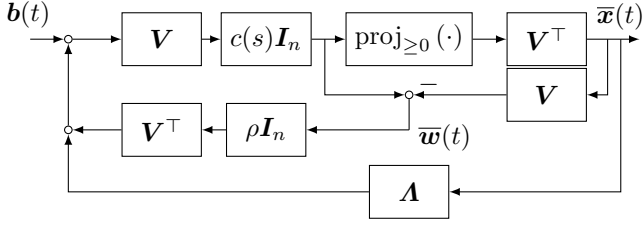
\begin{figure}[ht] 
\vskip 0.2in 
\begin{center} 
\centerline{\IfFileExists{./Figures/controller1.tex}{\input{./Figures/controllerTmpV2.tex}}{?}} 
\caption{Intermediate equivalent block diagram of Fig.~\ref{fig:alg0}, obtained by exploiting the eigendecomposition $\A = \V \LambdaB \VT$ and denoting by $\xbart$, $\bbart$, and $\wbart$ as in~\eqref{eq:transfSignals}.} 
\label{fig:alg1} 
\end{center} 
\vskip -0.2in 
\end{figure}

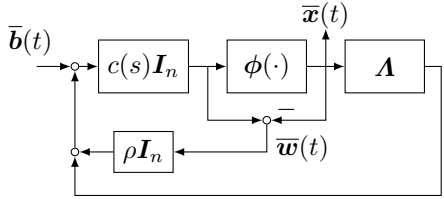
\begin{figure}[ht] 
\vskip 0.0in 
\begin{center} 
\centerline{\IfFileExists{./Figures/controllerDec.tex}{\input{./Figures/controllerDec.tex}}{?}} 
\caption{Final decoupled block diagram equivalent to Fig.~\ref{fig:alg0} and~\ref{fig:alg1} with the transformed signals $\xbart$, $\bbart$, and $\wbart$.}
\label{fig:algDec} 
\end{center} 
\vskip -0.2in 
\end{figure}

By leveraging the structure of the controller and the anti-windup mechanism, we can shift the matrix $\V$ after the controller block $\cs\In$, and $\VT$ before the anti-windup block $\rho\In$. Introducing the nonlinear mapping $\phiB:\Rn \to \Rn$, defined for all $\u \in \Rn$ as
\begin{equation}\label{eq:phiDef}
    \phiB{\left(\u\right)} \coloneqq \VT\projPlus{\left(\V\u\right)},
\end{equation}
and defining the transformed signals 
\begin{subequations}\label{eq:transfSignals}
    \begin{align}
        \xbart &\coloneqq \VT\xt,\label{eq:xbar}\\
        \bbart &\coloneqq \VT\bt,\label{eq:bbar}\\
        \wbart &\coloneqq \VT\wt,\label{eq:qbar}
    \end{align}
\end{subequations}
    we obtain the equivalent scheme shown in Fig.~\ref{fig:algDec}.

This structural simplification facilitates the subsequent analysis by enabling a decoupled representation of the system dynamics. Furthermore, it provides a principled framework to handle the uncertainty associated with the positive-definite matrix $\A$, allowing the control synthesis problem to be formulated within a robust Linear Matrix Inequality (LMI) setting.

\subsection{Step 1: \textit{Unconstrained} Convergence through the IMP}

\begin{figure}[ht]
    \vskip 0.0in
    \begin{center}
        \centerline{\IfFileExists{./Figures/controllerLin.tex}{\input{./Figures/controllerLin.tex}}{?}}
        \caption{Equivalent block diagrams of Fig.~\ref{fig:algDec} in the \textit{unconstrained} case, showing both vector and scalar representations.}
        \label{fig:algLin}
    \end{center}
    \vskip -0.2in
\end{figure}
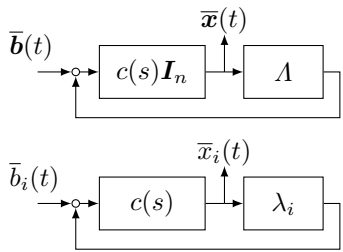
Following the approach in~\cite{bastianello_internal_2022}, the goal in the \textit{unconstrained} setting is to drive the gradient to zero, i.e., $\gradfxt \to 0$ as $t \to +\infty$. Due to the convexity of the cost function~\eqref{eq:quadratic-cost}, this condition implies that $\xt \to \xStart$, ensuring convergence to the optimal solution in the original closed-loop configuration of Fig.~\ref{fig:alg0}.

To be more precise, observe that
\begin{equation}\label{eq:uncoCase}
    \gradfxt = \A\xt + \bt = \V\left(\LambdaB\xbart + \bbart\right),
\end{equation}
where $\xbart$ and $\bbart$ are defined in~\eqref{eq:transfSignals}.
Ensuring that $\gradfxt \to 0$ as $t\to +\infty$ is equivalent to requiring
\begin{equation*}
    \LambdaB\xbart + \bbart \to 0,
\end{equation*}
which, due to the diagonal structure of $\LambdaB$, reduces to the scalar conditions
\begin{equation}\label{eq:errorto0}
    \ybarit \coloneqq \lambda_i \xbarit + \bbarit \to 0,
\end{equation}
for all $i = 1, \ldots, n$. This can be achieved by invoking the IMP. Indeed, analyzing the block diagram in Fig.~\ref{fig:algLin}, we see that the $i$-th decoupled component in the Laplace transform domain yields
\begin{equation}\label{eq:errorS}
    \Laplace\left[\ybarit\right] = \frac{1}{1 - \lambda_i \cs} \, \Laplace\left[\bbarit\right].
\end{equation}
By Assumption~\ref{as:modelbk}, the Laplace transform of $\bbarit$ admits a fraction representation as $\Laplace\left[\bbarit\right] = n_i(s) / \ds$, where $n_i(s)$ is a polynomial in $s$, and $\ds = \det(s \Imm - \F)$ is a common marginally stable denominator, shared by all $i = 1, \ldots, n$, inherited from the dynamics of $\bt$ given in \eqref{eq:linearModels}. By imposing $\cs = \ns/\ds$, for some suitable numerator $\ns$, the expression~\eqref{eq:errorS} becomes
\begin{equation}\label{eq:errorFinal}
    \Laplace\left[\ybarit\right] = \frac{\ds}{\ds - \lambda_i \ns} \cdot \frac{n_i(s)}{\ds}=\frac{n_i(s)}{\ds - \lambda_i \ns}.
\end{equation}
Hence, ensuring convergence $\ybarit \to 0$ as $t\to+\infty$ amounts to requiring that all poles of~\eqref{eq:errorFinal} have strictly negative real parts and this is guaranteed if we choose $\ns$ such that the polynomial $\ds - \lambda_i \ns$ is Hurwitz for all $\lm\leq\lambda_i \leq \lM$.
This robust stabilization condition can be addressed via several classical synthesis techniques. In this work, we adopt the continuous-time LMI-based procedure proposed in~\cite{luis_new_2015}, which generalizes the discrete-time framework of~\cite{de_oliveira_new_1999}.

Thus, the final controller that guarantees asymptotic tracking of the optimal trajectory $\xStart$ in the \textit{unconstrained} case—by virtue of the IMP~\cite{francis_internal_1976}—is given by $\Cs = \cs \In$ in Fig.~\ref{fig:algGen}, where
\begin{equation}\label{eq:scalarc}
    \cs = \K(s\In - \F)^{-1}\H,
\end{equation}
with $\F$ and $\H$ fixed by the internal model~\eqref{eq:linearModels}, and $\K$ synthesized via robust control design. Its corresponding full state-space realization reads:
\begin{equation}\label{eq:controllerEqua}
\begin{array}{rl}
    \detaextt &= \Fext \etaextt + \Hext \ut, \\
    \yct &= \Kext \etaextt,
\end{array}
\end{equation}
where $\Fext$ and $\Hext$ are defined in~\eqref{eq:linearModelsEst}, and $\Kext = \I_n \otimes \K$.

\subsection{Step 2: Anti-windup synthesis}\label{subsec:step2}
This subsection describes the second step of the synthesis procedure, which focuses on stabilizing the control scheme shown in Fig.~\ref{fig:algDec} in the \textit{constrained} case, i.e., when $\phiB{\left(\cdot\right)}\neq\In$.

The first step in this direction consists of categorizing the nonlinearity~\eqref{eq:phiDef} within the standard nonlinear framework introduced in~\cite{kothare_unified_1994}. To this end, we introduce the following definition from~\cite{grimm_linear_2004}.

\begin{definition}\label{def:PhiI}
    A function $\phiB:\Rn \to \Rn$ is said to \textit{belong to} $\Phi_\I$ if:
    \begin{enumerate}
        \item $\phiB(0) = 0$.\label{en:item3}
        \item $\phiB$ is locally Lipschitz continuous; \label{en:item1}
        \item $\phiB$ is such that
    \begin{equation}\label{eq:incrSec}
        \left(\jacphiu\right)^\top \left(\In - \jacphiu\right) \succeq 0
    \end{equation} for almost all $\u \in \Rn$, where $\jacphiu$ denotes the Jacobian of $\phiB$ evaluated at $\u$\footnote{For a non symmetric matrix $\A$ the symbol $\A\succeq 0$ means $\A+\A^T\succeq 0$.}.\label{en:item2}
        
    \end{enumerate}
\end{definition}
The next proposition establishes that the nonlinearity defined in~\eqref{eq:phiDef} satisfies Definition~\ref{def:PhiI}. This alignment ensures that our analysis is compatible with the hypotheses introduced in~\cite{grimm_linear_2004}, enabling us to apply the same theoretical framework developed therein.

\begin{prop}
    The function $\phiB{\left(\u\right)}$ defined in~\eqref{eq:phiDef} satisfies Definition~\ref{def:PhiI}.
\end{prop}

\begin{proof}
Condition~\ref{en:item3}) is trivially satisfied. Condition~\ref{en:item1}) is satisfied because $\projPlus{\left(\cdot\right)}$ corresponds to the element-wise application of the ReLU function~\cite{dubey_activation_2022}, which is Lipschitz continuous. Moreover, the composition of Lipschitz continuous functions is itself Lipschitz continuous~\cite{eriksson_applied_2010}.

For condition \ref{en:item2}), we compute the Jacobian of $\phiB$ at $\u$:
\begin{equation*}
\jacphiu = \VT \frac{\partial\projPlus{(\V\u)}}{\partial \u} = \VT \X \V,
\end{equation*}
where $\X = \mathrm{diag}(\indicatorPos(\v_1),\dots,\indicatorPos(\v_n))$ and $\v \coloneqq \V\u$, with $\indicatorPos(\cdot)$ denoting the nonnegative indicator function. Evaluating the left-hand side of the incremental sector condition~\eqref{eq:incrSec}:
\begin{align*}
\VT \X^\top \V(\In - \VT \X \V) &= \VT(\X^\top - \X^\top\X)\V \\ &= 0 \succeq 0,
\end{align*}
since $\X$ is diagonal with entries in $\{0,1\}$. Thus, condition \ref{en:item2}) holds, completing the proof.
\end{proof}

At this point, we can leverage the results of~\cite{grimm_linear_2004} to establish that the block scheme in Fig.~\ref{fig:algDec} guarantees a \textit{quadratic performance level} $\gamma$ for all $\phiB$ belonging to $\Phi_\I$. That is,
\begin{definition}\label{def:perfIndex}
    Given $\gamma > 0$, the system described in Fig.~\ref{fig:algDec} and Equation~\eqref{eq:controllerEqua} achieves a \textit{quadratic performance of level} $\gamma$ for all $\phiB$ belonging to $\Phi_\I$ if there exists a scalar $\epsilon > 0$ and a quadratic Lyapunov function $V(\etaext) = \etaext^\top \P \etaext$, such that the time derivative $\dot{V}(\etaext)$ along the dynamics described by Fig.~\ref{fig:algDec} and~\eqref{eq:controllerEqua} satisfies:
\begin{equation}
\label{eq:Vdot}
\begin{gathered}
    \dot{V}(\etaext) < -\epsilon\etaext^\top\etaext - \frac{1}{\gamma}\z^\top\z + \gamma\b^\top\b \\ \forall (\etaext,\b)\in\RNumb^{nm}\times\Rn,
\end{gathered} 
\end{equation}
where
\begin{equation}\label{eq:perfIndex}
    \z \coloneqq \gradf + \w,
\end{equation}
 quantifies the approximation accuracy\footnote{Requiring~\eqref{eq:Vdot} considering $\z$ and $\b$ or their decoupled counterpart $\zbar \coloneqq \VT\z = \LambdaB\xbar + \bbar + \wbar$ and $\bbar$ is equivalent due to the orthogonality of $\V$.} of $\xt$ to $\xStart$, and $\b$ represents the exogenous input~\eqref{eq:linearModelsEst}.
\end{definition}



\begin{rem}
    The choice of the performance index~\eqref{eq:perfIndex} in Definition~\ref{def:perfIndex} is motivated by the following optimization-theoretic interpretation. The Karush-Kuhn-Tucker (KKT) conditions, that are necessary and sufficient for the optimality in problem~\eqref{eq:problem1}, state that for any fixed $t \in \RNumb_+$, $\xStar(t) \in \Rn$ is the optimal solution iff there exists $\muStarB(t) \in \Rn$ such that 
    \begin{subequations}
    \begin{align}
        &\text{stationarity:}& &  \nabla f_t(\x^*(t)) - \muStarB(t) = 0,\label{eq:sta}\\
        &\text{primal feasibility:}& & \xStar(t)\geq 0,\label{eq:primFeas}\\
        &\text{dual feasibility:}& & \muStarB(t) \geq 0,\label{eq:dualFeas}\\
       &\text{complementary slackness:}& & \muStarB(t) \odot \xStar(t) = 0,\label{eq:compSl}
    \end{align}
    \end{subequations}
    
    These conditions hold in both the \textit{constrained} and \textit{unconstrained} settings. In particular, the complementary slackness condition~\eqref{eq:compSl}, together with~\eqref{eq:primFeas} and~\eqref{eq:dualFeas}, implies that for each component $i$, it must hold that either $\xStari = 0$ and $\muStari > 0$, or $\xStari > 0$ and $\muStari = 0$. Simultaneous positivity of both $\xStari$ and $\muStari$ is not admissible for any $i = 1,\dots,n$.

    Remarkably, the variables $-\wt$ and $\xt$, as defined in Fig.~\ref{fig:algGen}, inherently respect this structure. Specifically, for each $i$, the quantities $-\wit$ and $\xit$ are mutually exclusive: 
    \begin{itemize}
        \item If the constraint is inactive, we have \begin{equation*}
        \xit = \projPlus(\ycit) = \ycit > 0,
        \end{equation*}
        and thus $-\wit = \xit - \ycit = 0$;
        \item When the constraint is active then $\xit = 0$ (i.e., $\ycit < 0$), it follows that
        \begin{equation*}
            -\wit = \xit - \ycit > 0.
        \end{equation*}
    \end{itemize}
    Hence, the algorithm described in Fig.~\ref{fig:algGen} structurally enforces conditions~\eqref{eq:primFeas}–\eqref{eq:compSl} at all times, provided that $-\wt$ is interpreted as the vector of Lagrange multipliers associated with the nonnegativity constraint in problem~\eqref{eq:problem1}. 
     
 Under this interpretation, the performance output $\zt = \gradfxt - (-\wt)$ in~\eqref{eq:perfIndex} measures the deviation from the stationarity condition~\eqref{eq:sta}. Consequently, minimizing the norm of $\zt$ (or its integral over time) is a natural and meaningful performance objective.
\end{rem}

Verifying the inequality~\eqref{eq:Vdot} is sufficient to ensure stability and to certify a finite $\mathcal{L}_2$ gain $\gamma$ from the exogenous input $\b$ to the performance output $\z$, for all nonlinearities $\phiB$ belonging to $\Phi_\I$~\cite{grimm_linear_2004}.

We are now in a position to provide a formal result for computing the anti-windup gain $\rho$ via the following proposition.

\begin{prop}\label{prop:LMI}
        Given a desired performance level $\gamma > 0$, suppose there exist a matrix $\QB = \QBT \succ 0$, a scalar $\delta > 0$, and $\xi \neq 0$ such that the following LMIs hold for $\lambda = \lm$ and $ \lM$:
    {\small
    \begin{equation}
    \label{eq:LMI}
    \begin{split}
        &\hspace{-1.5mm}\left[\begin{array}{c@{\quad}c@{\quad}c@{\quad}c}\QB\Fcl^\top{+}\Fcl\QB & \H\left(\xi-\delta\lambda\right)+\QB\K^\top & \H & \QB\K^\top\lambda \\
        \star & -2\delta & 0 & \delta\left(1-\lambda\right)\\
        \star & \star &-\gamma &1\\
        \star &\star&\star & -\gamma
        \end{array}\right]\\
        &\qquad\prec 0
        \end{split}
    \end{equation}
    }
    where  $\Fcl \coloneqq \F + \H\lambda\K$. Then, the algorithm in Fig.~\ref{fig:algDec} with $\rho \coloneqq \xi/\delta$ satisfied Definition~\ref{def:perfIndex}. 
\end{prop}
\begin{proof}
    We first define $\P{\left(\lambda,\QB,\delta,\xi\right)}$ as the block matrix appearing in~\eqref{eq:LMI}, considered as a function of $\lambda$, $\QB$, $\delta$ and $\xi$. By assumption, we know that $\P{\left(\lm;\QB,\delta,\xi\right)} \prec 0$ and $\P{\left(\lM;\QB,\delta,\xi\right)} \prec 0$.

    It is well known that, for fixed $\QB$, $\delta$ and $\xi$, the set of $\lambda$ such that $\P{\left(\lambda,\QB,\delta,\xi\right)}\prec 0$ is convex~\cite{caverly_lmi_2019}. Consequently, $\P{\left(\lambda;\QB,\delta,\xi\right)} \prec 0$ holds for any $\lambda\in[\lm,\lM]$ and hence for all (unknown) values $\lambda = \lambda_i$, where $\lambda_i$ are the $n$ eigenvalues of $\LambdaB$.

    We are now in a position to apply ~\cite[Proposition 1]{grimm_linear_2004}. Indeed, the LMI appearing in~\cite[Proposition 1]{grimm_linear_2004}, when applied to the system described in Fig.~\ref{fig:algDec} with the controller given in ~\eqref{eq:controllerEqua}, admits as a solution the matrix $Q = \I_n \otimes \QB$, since this LMI is equivalent to the $n$ independent LMIs $\P{\left(\lambda_i;\QB,\delta,\xi\right)} \prec 0$ for each $i = 1, \ldots, n$.

    Since, by the previous argument, these LMIs are satisfied for all $\lm\leq\lambda_i \leq \lM$, the system in Fig.~\ref{fig:algDec}, with dynamics given by~\eqref{eq:controllerEqua}, satisfies the performance condition stated in Definition~\ref{def:perfIndex}, as guaranteed by~\cite[Proposition 1]{grimm_linear_2004}.
\end{proof}


Thus, the controller $\cs$ designed in \textit{Step 1}, together with the anti-windup gain $\rho$ computed via Proposition~\ref{prop:LMI}, ensures that the closed-loop algorithm in Fig.~\ref{fig:algDec} satisfies the performance criterion specified in Definition~\ref{def:perfIndex}, as desired.
\section{SIMULATIONS}

In this section, we compare the performance of three methods applied to problem~\eqref{eq:problem1}: the standard OP-GD algorithm~\cite{davydov_contracting_2023}, implemented as in~\eqref{eq:onlineProjDesc} with step size $\alpha = 1/\lM$; the continuous-time counterpart of the non-optimized \textit{Projected-Internal Model Anti-Windup} (non-optimized P-IMAW) proposed in~\cite{casti_control_2024}; and the P-IMAW introduced here in Fig.~\ref{fig:algGen}.

For both internal model-based algorithms, the robust controller is designed according to \textit{Step 1}, leveraging the result in~\cite[Theorem 4 with $\xi = 0.54$, $\epsilon = 1$]{luis_new_2015}. The parameter $\rho$ is computed via Proposition~\ref{prop:LMI} setting $\gamma = 10$.

In the first simulation, shown in Fig.~\ref{fig:trianWave}, we consider problem~\eqref{eq:problem1} with $n = 10$ and $\A = \V\LambdaB\VT$, where $\LambdaB$ is a diagonal matrix with eigenvalues randomly sampled from the uniform distribution over $\left[1,\,10\right]$, and $\V$ is a random orthogonal matrix. The signal $\bt$ is generated as a triangular wave~\cite[Fig. 3]{casti_control_2024} i.e., a piecewise linear ramp, with angular frequency $\omega = \qty[per-mode = symbol]{0.25}{\radian\per\second}$, so inducing several activations and deactivations of the constraints, allowing us to better appreciate the anti-windup behavior.

The internal model adopted for the design is associated with the polynomial
\begin{equation}\label{eq:linRamp}
    p(s) = s^2,
\end{equation}
and the pair $\left(\F,\,\H\right)$ is assumed in controllable canonical form with $\F$ having characteristic polynomial $\ps$, then corresponding to a linear ramp evolution.

The resulting tracking error $\norm{\xt - \xStart}$ is plotted over a $\qty[per-mode = symbol]{45}{\second}$ time window. As expected, the triangular wave introduces unavoidable transients due to constraint switching and only partial satisfaction of the IMP. Nonetheless, the improved performance of the proposed P-IMAW is evident, particularly during switching events generating windup effects.

\begin{figure}[ht] 
    \vskip 0.2in 
    \begin{center} 
    \centerline{\IfFileExists{./Figures/triangWave.tex}{\input{./Figures/triangWave.tex}}{?}} 
    \caption{Tracking error $\norm{\xt - \xStart}$ under triangular wave evolution of $\bt$ and internal model~\eqref{eq:linRamp}.} 
    \label{fig:trianWave} 
    \end{center} 
    \vskip -0.2in 
\end{figure}

In Fig.~\ref{fig:sineWave}, we report a second simulation under similar conditions. The objective function matrix $\A$ and the algorithm design is the same as in the previous experiment. This time, the signal $\bt$ is generated according to Assumption~\ref{as:modelbk}, with $\H$ and $\F$ as before but with 
\begin{equation}
    p(s) = \left(s^2 + \omega^2\right)s,
\end{equation}
which corresponds to a sinusoidal signal with angular frequency $\omega = \qty[per-mode = symbol]{0.25}{\radian\per\second}$ and an added constant component.

In this case, the improved tracking performance of the P-IMAW is even more evident, especially during transients triggered by constraint activation/deactivation.

\begin{figure}[ht] 
    \vskip 0.2in 
    \begin{center} 
    \centerline{\IfFileExists{./Figures/sineConst.tex}{\input{./Figures/sineConst.tex}}{?}} 
    \caption{Tracking error $\norm{\xt - \xStart}$ with $\bt$ following a sinusoidal plus constant evolution.} 
    \label{fig:sineWave} 
    \end{center} 
    \vskip -0.2in 
\end{figure}


\addtolength{\textheight}{-3cm}   


\section{CONCLUSIONS}
In this work, we proposed a novel method for online optimization with nonnegativite constraints. The method leverages the IMP and anti-windup techniques to improve the tracking performance in dynamic optimization settings. Future research directions include extending the proposed framework to more general scenarios, such as its integration with dual variable methods. Another promising direction is the adaptation of the same design principles to the direct discrete-time design setting.


\bibliographystyle{IEEEtran}
\bibliography{references}

\end{document}

%% file: Figures/tikzFigOptions.tex
\usetikzlibrary {chains,scopes,positioning,decorations.pathreplacing,calc,fit} 
\tikzset{token/.style = {draw=none, ultra thin, minimum size=5mm},
trigger/.style = {draw=none,fill = black,ultra thin, minimum size=5mm},
blank/.style = {draw, ultra thin, minimum size=5mm},
questionMark/.style = {draw, ultra thin, minimum size=5mm},
braced/.style={
     decoration={brace, mirror},
     decorate
   },
   braceu/.style={
     decoration=brace,
     decorate
   },
   dot/.style={
   	circle,inner sep=0pt,minimum size=0pt,fill=white
   	},
    block/.style={draw, fill=white, rectangle, 
            minimum height=2em, minimum width=3em},
    blockSmall/.style={draw, fill=white, rectangle, 
            minimum height=1.25em, minimum width=1.45em},
    input/.style={inner sep=0pt},       
    output/.style={inner sep=0pt},      
    sum/.style = {draw, fill=white, circle, minimum size = 0.1cm, inner sep=0pt},
    cross/.style={path picture={\draw[black] (path picture bounding box.south east) -- (path picture bounding box.north west) (path picture bounding box.south west) -- (path picture bounding box.north east);
}},
    mult/.style = {draw, fill=white, circle,cross, minimum size = 0.125cm, node distance=1.5cm, inner sep=0pt},
    pinstyle/.style = {pin edge={to-,thin,black}},
    sigmoid/.style = {blockSmall, path picture ={
    \draw[shift={(path picture bounding box.center)},
        domain=-0.3:0.3, samples=50,
    ] plot ({\x},{0.3*(sigma(30*\x)-0.5)});
}},
    declare function={
        sigma(\x) = 1/(1+exp(-\x));
    }
}

\definecolor{myBlue}{rgb}{0,0.4470,0.7410}
\definecolor{myRed}{rgb}{0.8500,0.3250, 0.0980}
\definecolor{myYellow}{rgb}{0.9290,0.6940,0.1250}
\definecolor{myPurple}{rgb}{0.4940, 0.1840, 0.5560}
\definecolor{myGreen}{rgb}{0.4660,0.6740,0.1880}

\definecolor{mycolor1}{rgb}{0.00000,0.44700,0.74100}%
\definecolor{mycolor2}{rgb}{0.85000,0.32500,0.09800}%
\definecolor{mycolor3}{rgb}{0.92900,0.69400,0.12500}%


%% file: notation.tex
\DeclareMathOperator*{\argmin}{arg\,min}

\DeclareMathOperator{\proj}{proj}

\def\projPlus{\proj_{\geq 0}}

\newcommand{\norm}[1]{\left\lVert#1\right\rVert}

\def\RNumb{{\mathbb{R}}}

\def\Rn{{\mathbb{R}^{n }}}
\def\Rm{{\mathbb{R}^{m }}}

\def\Rnn{{\mathbb{R}^{n\times n }}}
\def\C{{\mathbb{C}}}

\def\indicatorPos{\mathcal{X}_{\geq 0}}
\def\Laplace{\mathscr{L}}

\def\gradftdot{{\nabla f_t{\left(\cdot\right)}}}
\def\gradfqt{{\nabla \ftqt}}

\def\gradfxt{{\nabla \ftxt}}

\def\jacphiu{\frac{\partial \phiB}{\partial \u}{\left(\u\right)}}

\def\ftqt{f_t{\left(\qt\right)}}

\def\ftxt{f_t{\left(\xt\right)}}

\def\gradf{{\nabla \ftx}}
\def\ftx{f_t{\left(\x\right)}}

\def\xbart{\xbar{\left(t\right)}}
\def\xbarit{\overline{x}_i{\left(t\right)}}

\def\zbar{\overline{\z}}

\def\bbar{\overline{\b}}
\def\wbart{\wbar{\left(t\right)}}
\def\bbart{\bbar{\left(t\right)}}

\def\xbar{\overline{\x}}

\def\wbar{\overline{\w}}

\def\xint{\xiB_n{\left(t\right)}}
\def\xiOnet{\xiB_1{\left(t\right)}}
\def\xiit{\xiB_i{\left(t\right)}}
\def\dxiit{\dot{\xiB}_i{\left(t\right)}}

\def\xiextt{\xiB_{ext}{\left(t\right)}}
\def\dxiextt{\dot{\xiB}_{ext}{\left(t\right)}}

\def\qt{\qq{\left(t\right)}}
\def\ft{\f{\left(t\right)}}
\def\et{\e{\left(t\right)}}

\def\wt{\w{\left(t\right)}}
\def\vt{\v{\left(t\right)}}
\def\ut{\u{\left(t\right)}}
\def\uct{\u_c{\left(t\right)}}
\def\yct{\y_c{\left(t\right)}}

\def\etaext{\etaB_{ext}}
\def\etaextt{\etaB_{ext}{\left(t\right)}}

\def\detaextt{\dot{\etaB}_{ext}{\left(t\right)}}

\def\Cs{\C{\left(s\right)}}
\def\Cxps{\C_{\x\p}{\left(s\right)}}
\def\Cxqs{\C_{\x\qq}{\left(s\right)}}
\def\Cyps{\C_{\y\p}{\left(s\right)}}
\def\Cyqs{\C_{\y\qq}{\left(s\right)}}

\def\ps{p{\left(s\right)}}
\def\ds{d{\left(s\right)}}
\def\ns{n{\left(s\right)}}
\def\cs{c{\left(s\right)}}

\def\yt{\y{\left(t\right)}}
\def\st{\s{\left(t\right)}}
\def\xt{\x{\left(t\right)}}
\def\pt{\p{\left(t\right)}}
\def\qt{\qq{\left(t\right)}}
\def\dqt{\dot{\qq}{\left(t\right)}}

\def\VT{\V^\top}

\def\Fext{\F_{ext}}

\def\Hext{\H_{ext}}

\def\Kext{\K_{ext}}

\def\xStari{x^*_i}
\def\xStar{\x^*}
\def\zt{\z{\left(t\right)}}
\def\bt{\b{\left(t\right)}}

\def\bTt{\b^\top{\left(t\right)}}
\def\xStart{\xStar{\left(t\right)}}
\def\ct{c{\left(t\right)}}

\def\ybarit{\overline{y}_i{\left(t\right)}}

\def\bbarit{\overline{b}_i{\left(t\right)}}
\def\bit{b_i{\left(t\right)}}

\def\wit{w_i{\left(t\right)}}
\def\xit{x_i{\left(t\right)}}
\def\ycit{{y_c}_i{\left(t\right)}}

\def\QB{\overline{\Q}}
\def\QBT{\QB^\top}
\def\lm{{\lambda_{\min}}}
\def\lM{{\lambda_{\max}}}

\def\Fcl{\F{\left(\lambda\right)}}

\def\In{\I_n}
\def\Imm{\I_m}

\def\b{{\mathbold{b}}}

\def\e{{\mathbold{e}}}
\def\f{{\mathbold{f}}}

\def\p{{\mathbold{p}}}
\def\qq{{\mathbold{q}}}
\def\s{{\mathbold{s}}}
\def\v{{\mathbold{v}}}
\def\w{{\mathbold{w}}}
\def\v{{\mathbold{v}}}
\def\x{{\mathbold{x}}}
\def\y{{\mathbold{y}}}
\def\z{{\mathbold{z}}}
\def\u{{\mathbold{u}}}

\def\A{{\mathbold{A}}}
\def\B{{\mathbold{B}}}
\def\C{{\mathbold{C}}}

\def\F{{\mathbold{F}}}

\def\H{{\mathbold{H}}}
\def\I{{\mathbold{I}}}

\def\P{{\mathbold{P}}}
\def\Q{{\mathbold{Q}}}

\def\V{{\mathbold{V}}}

\def\X{{\mathbold{X}}}

\def\K{{\mathbold{K}}}

\def\TT{{\mathbold{T}}}

\def\xiB{{\mathbold{\xi}}}
\def\etaB{{\mathbold{\eta}}}
\def\phiB{{\mathbold{\phi}}}

\def\muB{{\mathbold{\mu}}}

\def\muStarB{\muB^*}
\def\muStari{\mu^*_i}
\def\LambdaB{{\mathbold{\Lambda}}}

%% file: Figures/controllerGeneral.tex
\begin{tikzpicture}[auto, >=latex]
\node[input] (w){};
\node[input, above = 0.15cm of w] (wDelta){};
\node[input, below = 0.15cm of w] (u){};
\node [block, minimum width = 1.75cm,minimum height = 1 cm, right = of w] (P) {$\begin{bmatrix}
   \Cxps & \Cyps\\
   \Cxqs & \Cyqs
\end{bmatrix}$};
\node [block, above = 0.25cm of P] (Delta) {$\projPlus{\left(\cdot\right)}$};
\node [block, below = 0.25cm  of P] (ctrl) {$\gradftdot$};
\node[output,right = of P] (e){};
\node[output, above = 0.15cm of e] (eDelta){};
\node[output, below = 0.15cm of e] (y){};
\coordinate (tmpInputDelta) at ($(P.west)!0.5!(P.north west)$);
\coordinate (tmpInputU) at ($(P.west)!0.5!(P.south west)$);
\coordinate (tmpMiddleInput) at ($(w)!0.5!(P.west)$);
\coordinate (tmpMiddleOutput) at ($(e)!0.5!(P.east)$);
\coordinate (tmpOutputDelta) at ($(P.east)!0.5!(P.north east)$);
\coordinate (tmpOutputU) at ($(P.east)!0.5!(P.south east)$);

\coordinate (tmpMiddleInputUp) at (tmpInputDelta-|tmpMiddleInput);
\coordinate (tmpMiddleOutputUp) at (tmpOutputDelta-|tmpMiddleOutput);

\coordinate (tmpMiddleInputDown) at (tmpInputU-|tmpMiddleInput);
\coordinate (tmpMiddleOutputDown) at (tmpOutputU-|tmpMiddleOutput);

\draw [draw,->] (tmpOutputDelta) -- (tmpMiddleOutputUp)  |- node[below right]{$\pt$} (Delta);
\draw [draw,->] (tmpOutputU) -- (tmpMiddleOutputDown) |- node[above right]{$\qt$} (ctrl);
\draw [draw,->] (Delta) -| node[below left]{$\xt$} (tmpMiddleInputUp) -- (tmpInputDelta) ;
\draw [draw,->] (ctrl)-|node[above left]{$\yt$}  (tmpMiddleInputDown) -- (tmpInputU);
\end{tikzpicture}

%% file: Figures/controllerMainSheme.tex
\begin{tikzpicture}[ every join/.style = {-latex},node distance=1.25mm and .75cm]
{[start chain=trunk going below]
    \node[block,on chain] (ctrl){$\Cs$};
    {[start branch = tmpL going left]
        \coordinate[on chain,join = by latex-, label = {above:$\ut$}] (tmpL);
    }
    \node[block,on chain = going right,join = by {to path = {(\tikztostart) -- node[above]{$\pt$}(\tikztotarget)}}] (proj) {$\projPlus{\left(\cdot\right)}$};
    \coordinate (tmpProjL) at ($(proj.west) - (0.375,0)$) ;
    \coordinate (tmpProjR) at ($(proj.east) + (0.375,0)$) ;
    {[start branch = tmpR going right]
        \coordinate[on chain,join = by -] (tmpR);
        \node[output,on chain,join] {$\xt$};
    }
    \node[sum,on chain,yshift = -1mm,join = with tmpProjL by {to path = {(\tikztostart) |- (\tikztotarget)}},join = with tmpProjR by {to path = {(\tikztostart) |- (\tikztotarget)node[above right, yshift = -0.1cm]{$-$}}}] (sumAnt) {};

    \node[block,yshift = -1.2cm,on chain,join = with tmpR by {to path = {(\tikztostart) |- (\tikztotarget)}}] (oracle) {$\gradftdot$};
    \coordinate (tmpSumAnt) at ($(sumAnt)!0.5!(oracle.north)$);
    \node[sum] at ($(tmpL |- tmpSumAnt)$) (sumUc){};
    \node[block] at (sumUc -| ctrl) (rho){$\TT$};
    \draw [-latex] (sumAnt) -- node[below right]{$\wt$} (tmpSumAnt) -- (rho);
    \draw [-latex] (rho)  -- node[above]{$\vt$} (sumUc);
    \draw [-latex] (oracle) -| node[ above right]{$\yt$}(sumUc); 
    \draw [-] (tmpL) -- (sumUc); 
    \node[fit=(rho)(tmpL)(tmpR)(ctrl),inner xsep=0.5cm,inner ysep=0.2cm,gray, draw, dashed,  label = {[gray]north:P-IMAW algorithm}] {};
}
\end{tikzpicture}

%% file: Figures/controllerTmp.tex
\begin{tikzpicture}[ every join/.style = {-latex},node distance=1.25mm and .75cm]
{[start chain=trunk going below]
    \node[block,on chain] (ctrl){$\cs\In$};
    {[start branch = tmpL going left]
        \node[sum,on chain,join = by {to path = { [latex-] (\tikztostart)  -- node[above,xshift = -0.1cm]{$\ut$}(\tikztotarget)}} ] (tmpL) {};
        \node[input,on chain ,join = by latex-, label = {above:$\bt$}]{};
    }
    \node[block,on chain = going right,join = by {to path = {(\tikztostart) -- node[above]{$\pt$}(\tikztotarget)}}] (proj) {$\projPlus{\left(\cdot\right)}$};
    \coordinate (tmpProjL) at ($(proj.west) - (0.375,0)$) ;
    \coordinate (tmpProjR) at ($(proj.east) + (0.375,0)$) ;
    {[start branch = tmpR going right]
        \coordinate[on chain,join = by -] (tmpR);
        \node[output,on chain,join] {$\xt$};
    }
    \node[sum,on chain,yshift = -1mm,join = with tmpProjL by {to path = {(\tikztostart) |- (\tikztotarget)}},join = with tmpProjR by {to path = {(\tikztostart) |- (\tikztotarget)node[above right, yshift = -0.1cm]{$-$}}}] (sumAnt) {};

    \node[block,yshift = -0.75cm,on chain,join = with tmpR by {to path = {(\tikztostart) |- (\tikztotarget)}}] (oracle) {$\A$};
    \coordinate (tmpSumAnt) at ($(sumAnt)!0.5!(oracle.north)$);
    \node[sum] at ($(tmpL |- tmpSumAnt)$) (sumUc){};
    \node[block] at (sumUc -| ctrl) (rho){$\rho\In$};
    \draw [-latex] (sumAnt) -- node[below right]{$\wt$} (tmpSumAnt) -- (rho);
    \draw [-latex] (rho) -- node[above]{$\vt$} (sumUc);
    \draw [-latex] (oracle) -| node[ above right]{$\yt$} (sumUc);
    \draw [-latex] (sumUc) -- (tmpL);

}
\end{tikzpicture}

%% file: Figures/controllerTmpV2.tex
\begin{tikzpicture}[ every join/.style = {-latex},node distance=1.25mm and .45cm]
{[start chain=trunk going below]
    \node[block,on chain] (ctrl){$\cs\In$};
    {[start branch = tmpL going left]
        \node[block,on chain,join = by latex-,xshift = 0.2cm] (V1) {$\V$};
        \node[sum,on chain,xshift = -0.2cm,join = by {to path = { [latex-] (\tikztostart)  -- node[above,xshift = -0.1cm]{}(\tikztotarget)}} ] (tmpL) {};
        \node[input,on chain ,join = by latex-, label = {above:$\bbart$}]{};
    }
    \node[block,on chain = going right,join] (proj) {$\projPlus{\left(\cdot\right)}$};

    \coordinate (tmpProjL) at ($(proj.west) - (0.275,0)$) ;
   
    {[start branch = tmpR going right]
        \node[block,on chain = going right,join] (VT) {$\VT$};
        \coordinate[on chain,join = by -] (tmpR);
        \node[output,on chain,join = by {to path = {(\tikztostart) node[above ]{$\xbart$} -- (\tikztotarget)}},xshift = -0.2cm] {};
         \coordinate (tmpProjR) at ($(VT.east) + (0.275,0)$) ;
    }
    \node[sum,on chain,yshift = -2mm,join = with tmpProjL by {to path = {(\tikztostart) |- (\tikztotarget)}}] (sumAnt) {};
    
    {[start branch = branchVTBelow placed {at = (sumAnt -| VT)} ]
        \node[block, on chain, join = with sumAnt by {to path = { (\tikztotarget)  -- (\tikztostart)node[above right, yshift = -0.1cm]{$-$}}}, join = with tmpProjR by {to path = {(\tikztostart) |- (\tikztotarget)}}] (V){$\V$};
    }
    
    \node[block,yshift = -0.75cm,on chain,join = with tmpR by {to path = {(\tikztostart) |- (\tikztotarget)}}] (oracle) {$\LambdaB$};
    \coordinate (tmpSumAnt) at ($(sumAnt)!0.5!(oracle.north)$);
    \node[sum] at ($(tmpL |- tmpSumAnt)$) (sumUc){};
    \node[block] at (sumUc -| ctrl) (rho){$\rho\In$};
    \node[block] at (rho -| V1) (VT2){$\VT$};
    \draw [-] (sumAnt) -- node[below right]{$\wbart$}  (tmpSumAnt);
    \draw [-latex] (tmpSumAnt) -- (rho);
    \draw [-latex] (rho) -- (VT2) ;
    \draw [-latex] (VT2) -- node[above]{} (sumUc);
    \draw [-latex] (oracle) -| (sumUc);
    \draw [-latex] (sumUc) -- (tmpL);
}
\end{tikzpicture}

%% file: Figures/controllerDec.tex
\begin{tikzpicture}[ every join/.style = {-latex},node distance=5mm and .25cm]
{[start chain=trunk]
    \node [input,on chain,label = {[shift={(-0.5,0)}] above right:$\bbart$ }] (input) {};
    \node[sum,on chain,xshift = 0.2cm,join ] (sumGrad){};
    \node[block,on chain,join= by {to path={(\tikztostart) -- (\tikztotarget) }}] (sysController) {$\cs\In$};
    {[start branch = rhoBranch going below]
        \node[blockSmall, on chain] (rho) {$\rho\In$};
        {[start branch = rhoBelow going below ]
            \coordinate[on chain,yshift = 2mm] (rhoBelow);
        }
    }
    {[start branch = sumBranch placed {at = (sumGrad |- rho)} ]
        \node[sum, on chain, join = with rho, join =  with sumGrad by latex-] (sumRho){};
    }
    
    \coordinate[on chain , join = by -] (yc);
    \node[block,on chain,join=by {to path={(\tikztostart) node[above]{} --  (\tikztotarget) }}] (project) {$\phiB{\left(\cdot\right)}$};
    {[start branch = antiWindupBranch going below ]
        \node[sum,on chain,yshift = +.2cm,join = with yc by {to path ={(\tikztostart) |- (\tikztotarget)}}] (q){};
        \chainin (rho)[join = by {to path = {(\tikztostart) |- node[below right, yshift = 10]{$\wbart$}(\tikztotarget)}}] ;
    }
    \coordinate[on chain = trunk, join = by -] (x);
    {[start branch = xout going above]
    \node[output,on chain,join]{$\xbart$};
    }
    {[continue branch = antiWindupBranch]
        \chainin (q) [join = with x by {to path ={(\tikztostart) |- (\tikztotarget) node[above right, yshift = -0.1cm]{$-$}}}];
    }
    \node[block,on chain = trunk,join = by {to path={(\tikztostart)--  (\tikztotarget) }}] (plant) {$\LambdaB$};
    \chainin (sumRho) [join = by {to path = {(\tikztostart) -- node[above right, yshift = 0]{}($(\tikztostart) + (0.75,0)$) |- (rhoBelow) -|(\tikztotarget)}}];
}
\end{tikzpicture}

%% file: Figures/controllerLin.tex
\begin{tikzpicture}[ every join/.style = {-latex},node distance=2.5mm and .25cm]
{[start chain=trunk,local bounding box=scope1, anchor=center,align = center]
    \node [input,on chain,label = {[shift={(-0.5,0)}]above right:$\bbart$ }] (input) {};
    \node[sum,on chain,xshift = 0.2cm,join ] (sumGrad){};
    \node[block,minimum width=4em,on chain,join= by {to path={(\tikztostart) -- (\tikztotarget) }}] (sysController) {$\cs\In$};
    {[start branch = cBelow going below ]
            \coordinate[on chain] (cBelow);
    }
    \coordinate[on chain = trunk,join = by {to path={[-](\tikztostart)--  (\tikztotarget) }}] (x);
    {[start branch = out going above]
    \node[output,on chain,join,above,yshift = 0.25cm] {$\xbart$};
    }
    \node[block,on chain = trunk,join = by {to path={(\tikztostart)--  (\tikztotarget) }}] (plant) {$\Lambda$};
    \chainin (sumGrad) [join = by {to path = {(\tikztostart) -- ($(\tikztostart) + (0.75,0)$) |- (cBelow) -|(\tikztotarget)}}];
}
{[start chain=trunk1,local bounding box=scope2, shift={(0,-1.75cm)},anchor = center,align = center]
    \node [input,on chain=trunk1,label = {[shift={(-0.5,0)}]above right:$\bbarit$ }] (input1) {};
    \node[sum,on chain=trunk1,xshift = 0.2cm,join ] (sumGrad1){};
    \node[block,minimum width=4em ,on chain=trunk1,join= by {to path={(\tikztostart) -- (\tikztotarget) }}] (sysController1) {$\cs$};
    {[start branch = cBelow1 going below ]
            \coordinate[on chain] (cBelow1);
    }
    \coordinate[on chain = trunk1,join = by {to path={[-](\tikztostart)--  (\tikztotarget) }}] (x1);
    {[start branch = out1 going above]
    \node[output,on chain,join,above,yshift = 0.25cm]{$\xbarit$};
    }
    \node[block,on chain ,join = by {to path={(\tikztostart)--  (\tikztotarget) }}] (plant1) {$\lambda_i$};
    \chainin (sumGrad1) [join = by {to path = {(\tikztostart) -- ($(\tikztostart) + (0.75,0)$) |- (cBelow1) -|(\tikztotarget)}}];
}
\end{tikzpicture}